\documentclass[12pt,reqno]{amsart}

\usepackage{comment}
\usepackage{amsmath,amsthm,amscd,amsfonts,amssymb,mathtools}
\usepackage{mathrsfs}
\usepackage{graphicx,xcolor}
\usepackage{enumitem}
\usepackage[bookmarksnumbered,colorlinks,plainpages]{hyperref}

\newcommand{\N}{\mathbb N}
\newcommand{\R}{\mathbb R}
\newcommand{\C}{\mathbb C}

\renewcommand{\d}{\,\mathrm{d}}

\newcommand\haar{\mathbf{h}}

\newcommand\IAA{\mathbf{IA}_0}

\newtheorem{theorem}{Theorem}[section]
\newtheorem{lemma}[theorem]{Lemma}

\newtheorem{corollary}[theorem]{Corollary}
\theoremstyle{definition}

\theoremstyle{remark}
\newtheorem{remark}[theorem]{Remark}
\numberwithin{equation}{section}

\begin{document}
\setcounter{page}{1}

\color{black}

\title[Khintchine's limit theorem for monotone convolutions]{Khintchine's limit theorem for monotone convolutions}

\author[Cheng-Xuan Soong]{Cheng-Xuan Soong}

\address{
National Center for Theoretical Sciences, Mathematics Division, Taipei 106, Taiwan
}

\email{cxsoong@ncts.ntu.edu.tw}

\author[Takahiro Hasebe]{Takahiro Hasebe}

\address{
Department of Mathematics, Hokkaido University
Kita 10, Nishi 8, Kita-Ku, Sapporo, Hokkaido, 060-0810, Japan
}

\email{thasebe@math.sci.hokudai.ac.jp}

\date{\today}

\begin{abstract}
We prove that every weak limit of monotone convolution powers
is monotonically infinitely divisible. In contrast, we prove that every
probability measure on $\R$ can occur as a weak limit of monotone
convolutions of an infinitesimal triangular array. We also discuss similar limit theorems for multiplicative monotone convolution on the unit circle. 


\end{abstract}
\maketitle

\section{Introduction}

Our main purpose is to show analogues of Khintchine's limit theorem for additive monotone convolution $\rhd$ and multiplicative monotone convolution $\circlearrowright$. Let $\mathbf P(\R)$ denote the set of probability measures on $\R$.  For $\mu_n,\mu \in \mathbf{P}(\R),~n\in\N$, the weak convergence of $\mu_n$ to $\mu$ is denoted as $\mu_n\Rightarrow \mu$.  We denote by $\mathbf{IA}(\rhd)$ the class of weak limits arising from the additive monotone convolution of  
infinitesimal triangular arrays, i.e., 
$\mathbf{IA}(\rhd)$ is the collection of all $\mu \in \mathbf{P}(\R)$ such that there exists a family $\{\mu_{n,j}: 1\le j \le k_n, ~n\in\N\}$ in $\mathbf{P}(\R)$ with $k_n\in\N$, $k_n\to\infty$ and 
\begin{align*}
&\lim_{n\to\infty}\max_{1\le j \le k_n}\mu_{n,j}(\R\setminus[-\epsilon,\epsilon]) =0, \qquad \forall \epsilon>0, \qquad \text{(infinitesimality)}\\
&\mu_{n,1}\rhd \mu_{n,2}\rhd \cdots \rhd \mu_{n,k_n} \Rightarrow \mu. 
\end{align*}
We also denote by $\IAA(\rhd)$ the
class of weak limits arising from the monotone convolution of triangular arrays with identical entries in each row, i.e., $\IAA(\rhd)$ is the collection of all $\mu \in \mathbf{P}(\R)$ such that there exist  sequences $\{\mu_n\}_{n\in\N}$ in $\mathbf{P}(\R)$ and $\{k_n\}_{n\in\N}$ in $\N$ with $k_n\to\infty$ and 
\[
\mu_n^{\rhd k_n} \Rightarrow \mu. 
\]
Note that we do not assume the infinitesimality $\mu_n\Rightarrow \delta_0$ a priori, but this automatically follows as we shall see in Step 1 of Section \ref{sec:proof} in the proof of Theorem \ref{thm:main} below. Consequently, the class  $\IAA(\rhd)$ is a subclass of $\mathbf{IA}(\rhd)$. 

Let $\mathbf{ID}(\rhd)$ denote the class of
$\rhd$-infinitely divisible probability measures on $\R$.
The very definition of  infinite
divisibility implies 
\[
    \mathbf{ID}(\rhd)\subset \IAA(\rhd). 
\]
It is also known that
\[
    \mathbf{ID}(\rhd)\subsetneq \mathbf{Univ}(\R)
    \subset \mathbf{IA}(\rhd)\subset \mathbf P(\R), 
\]
where 
$\mathbf{Univ}(\R)$ denotes the class of probability measures whose
reciprocal Cauchy transforms are univalent on $\C^+$. 
The inclusion
$\mathbf{ID}(\rhd)\subset\mathbf{Univ}(\R)$ follows from the
semigroup characterization and the theory of holomorphic semigroups;
see \cite{Mur00,Bel05,BP78}, and its strictness follows from
\cite[Examples~6.12--6.14]{FHS20}. The inclusion
$\mathbf{Univ}(\R)\subset\mathbf{IA}(\rhd)$ follows from Loewner theory for holomorphic dynamics   
\cite[Theorem~6.9(1)]{FHS20}, while the last inclusion is immediate from
the definition.

For classical and free convolutions, the weak limits of convolutions of infinitesimal
triangular arrays are precisely the infinitely divisible probability
measures; see \cite{GK54,BP00}. Limit theorems for monotone convolution
have also been studied in several directions; see, for example,
\cite{Mur01,AW14,Wan12,Wan14}. Motivated by these classical, free and monotone 
results, Franz, Hasebe, and Schlei{\ss}inger conjectured that the classes 
$\IAA(\rhd)$ and $\mathbf{ID}(\rhd)$ coincide; 
see \cite[Conjecture~6.39]{FHS20}. We prove this conjecture.

\begin{theorem}
\label{thm:main}
We have
\[
    \IAA(\rhd)=\mathbf{ID}(\rhd).
\]
\end{theorem}

Applying Theorem~\ref{thm:main} to $n$th monotone convolution
roots of a weakly convergent sequence in $\mathbf{ID}(\rhd)$
gives the following corollary, answering
\cite[Conjecture~6.17]{FHS20} affirmatively.

\begin{corollary}\label{cor:weak-closedness}
The class $\mathbf{ID}(\rhd)$ is closed under weak convergence.
\end{corollary}

In contrast, weak limits of the monotone convolutions of general infinitesimal  triangular
arrays need not be monotonically  infinitely divisible; see
\cite[Theorem~6.9(1) and Example~6.12]{FHS20}.
Franz, Hasebe, and Schlei{\ss}inger further conjectured that $
    \mathbf{IA}(\rhd)=\mathbf{Univ}(\R)$; see \cite[Conjecture~6.10]{FHS20}. This conjecture was supported by
several partial results. In particular, if the maximal variance in each
row tends to zero, then every weak limit has a univalent reciprocal
Cauchy transform; see \cite[Theorem~6.9(2)]{FHS20}. A partial converse
also holds for infinitesimal triangular arrays whose entries belong to
the class $\mathcal L$ introduced by Anshelevich and Arizmendi; see
\cite[Remark~51]{AA17} and \cite[Remark~6.11]{FHS20}. Moreover, a
closely related characterization holds on the unit circle; see
\cite[Theorem~7.5]{FHS20}.

Nevertheless, this conjecture is false. We prove the following.

\begin{theorem}
\label{thm:ia}
We have
\[
    \mathbf{IA}(\rhd)=\mathbf P(\R).
\]
\end{theorem}

In other words, weak limits  of monotone convolutions of infinitesimal triangular arrays impose no restriction
on the limit law. In this sense, the situation is analogous to the
Boolean case, where every probability measure is infinitely divisible;
see \cite{SW97}.
For instance, the symmetric Bernoulli measure gives an explicit element of $\mathbf{IA}(\rhd)\setminus\mathbf{Univ}(\R)$; see
Remark~\ref{rem:exp}, where we also explain why the corresponding
construction does not contradict the partial results above.

Combining the two theorems with the known inclusions, we obtain
\begin{equation}\label{eq:inclusions}
     \IAA(\rhd)
    =\mathbf{ID}(\rhd)
    \subsetneq\mathbf{Univ}(\R)
    \subsetneq\mathbf{IA}(\rhd)
    =\mathbf P(\R).
\end{equation}

We then study similar limit theorems for multiplicative monotone convolution on the unit circle. 
The analogue of Theorem \ref{thm:main} turns out to hold for multiplicative monotone convolution. On the other hand, the analogue of Theorem \ref{thm:ia} fails to hold according to \cite[Theorem~7.5]{FHS20}. For further details, see Theorem  \ref{thm:unit} and Remark \ref{rem:unit}. 

The structure of the paper is what follows. Section \ref{sec:pre} collects  definitions of basic notions and tools concerning monotone convolution $\rhd$. Sections \ref{sec:proof} and \ref{sec:proof2} offer proofs of Theorems \ref{thm:main} and \ref{thm:ia}, respectively. Section \ref{sec:unit} is devoted to similar limit theorems for multiplicative monotone convolution. 

\subsection*{Acknowledgement} This work was carried out during the authors' stays at National Tsing Hua University. The authors would like to thank Hao-Wei Huang for his hospitality.

\section{Preliminaries}\label{sec:pre}
Recall that the additive monotone convolution $\mu\rhd\nu$ is the distribution of the sum $X+Y$ where $X, Y$ are monotonically independent self-adjoint operators having distributions $\mu, \nu$, respectively \cite{Fra09}. We refer to \cite{Fra09,FHS20,Has,Mur00} for background on monotone convolution.

For $\lambda\in\mathbf{P}(\R)$, define its Cauchy transform and reciprocal Cauchy transform by
\[
G_\lambda(z):=\int_\R\frac{1}{z-s}\d\lambda(s),
\quad
F_\lambda(z):=\frac{1}{G_\lambda(z)},
\quad z\in\C^+.
\]
For $\mu,\nu\in\mathbf{P}(\R)$, their monotone convolution $\mu\rhd\nu$ is the unique probability measure characterized by
\[
F_{\mu\rhd\nu}=F_\mu\circ F_\nu.
\]
In particular, $F_{\mu^{\rhd n}}=F_\mu^{\circ n}$, where $F_\mu^{\circ n}$ denotes the $n$-fold composition of $F_\mu$. 

A measure $\mu\in\mathbf P(\R)$ is called $\rhd$-infinitely divisible if, for every $n\in\N$, there exists $\mu_n\in\mathbf P(\R)$ such that
\[
\mu=\mu_n^{\rhd n}.
\]
The $n$th roots $\mu_n$ are unique \cite[Proposition 3.8]{Bel05}. The class of all $\rhd$-infinitely divisible measures is denoted by $\mathbf{ID}(\rhd)$. 


We recall an analytic characterization of reciprocal Cauchy transforms. A function $f \colon \C^+\to \C$ is said to have a non-tangential limit $w \in \C \cup\{\infty\}$  at $\infty$ if for any $\alpha>0$,  
\[
\lim_{\substack{z\to\infty\\ \alpha|\Re z|<\Im z}} f(z) = w. 
\]
In this case, we write 
\[
w= \angle\!\lim_{z\to\infty} f(z).
\]
An analytic function $F\colon \C^+\to \C^+\cup \R$ has the Nevanlinna representation \begin{equation}
F(z)
=
a + b z + \int_\R\frac{1+sz}{s-z}\d\rho(s),
\qquad z\in\C^+ 
\label{eq:nevan}
\end{equation}
for some $a \in \R$, $b\ge0$  and finite positive measure $\rho$ on $\R$. From this formula, one can verify that the limits 
\[
\angle\!\lim_{z\to\infty} \frac{F(z)}{z}, \quad \angle\!\lim_{z\to\infty} F'(z), \qquad \lim_{y\to+\infty} \frac{F(iy)}{iy}
\]
all exist in $[0,+\infty)$ and coincide with $b$. We will write the common value as $b=F'(\infty)$.

A characterization of the reciprocal Cauchy transforms is given as follows \cite[Proposition 5.2]{BV93}. 
\begin{lemma}\label{lem:F}
    Let $F:\C^+\to\C^+\cup\R$ be analytic. The following are equivalent. 
    \begin{enumerate}[label=\rm(\arabic*)]
        \item\label{itemF1} $F=F_\mu$ for some $\mu \in \mathbf{P}(\R)$. 
        \item\label{itemF2} $F'(\infty)=1$. 
    \end{enumerate}
    Moreover, if \ref{itemF1} and \ref{itemF2} hold then we have $\Im F(z) \ge \Im z$ on $\C^+$.  
\end{lemma}


We will also use the standard relation between weak convergence and convergence of reciprocal Cauchy transforms; see \cite[Theorem~2.5]{Maa92} or \cite[Lemma~2.6]{FHS20}.

\begin{lemma}\label{lem:Fconv}
    Let $S$ be a subset of $\C^+$ having an accumulation point in $\C^+$. For $\lambda_n,\lambda\in\mathbf P(\R)$~$(n\in\N)$, the following are equivalent. 
    \begin{enumerate}[label=\rm(\arabic*)]
    \item $ \lambda_n\Rightarrow\lambda$, 
    \item $ F_{\lambda_n}\longrightarrow F_\lambda$ locally uniformly on
$\C^+$, 
\item $F_{\lambda_n}\longrightarrow F_\lambda$ pointwise on $S$. 
    \end{enumerate}
    \end{lemma}

Tightness can also be characterized by the reciprocal Cauchy transforms, see \cite[Proposition 4.39]{Has}.    \begin{lemma}\label{lem:Ftight}
    For a family $\mathcal{P}$ of probability measures, the following are equivalent. 
    \begin{enumerate}[label=\rm(\arabic*),itemsep=1mm]
    \item $ \mathcal{P}$ is tight. 
    \item $\displaystyle\lim_{y\to+\infty}\sup_{\mu \in \mathcal{P}}\left|\frac{F_\mu(iy)}{iy}-1\right| =0$.   
    \end{enumerate}
    \end{lemma}


We finally record the Schwarz--Pick inequality on the upper half-plane; see \cite[p.~136, Exercise~2]{Ahl79}.

\begin{lemma}[Schwarz--Pick lemma]
\label{lem:schwarz}
Let $f:\C^+\to\C^+$ be analytic. Then, for $z,w\in\C^+$,
\[
\frac{|f(z)-f(w)|}{|f(z)-\overline{f(w)}|}
\le
\frac{|z-w|}{|z-\overline{w}|}.
\]
Equivalently,
\[
\frac{|f(z)-f(w)|^2}{\Im f(z)\Im f(w)}
\le
\frac{|z-w|^2}{\Im z\Im w}.
\]
\end{lemma}

\section{Infinite divisibility of iid limits}\label{sec:proof}
We prove Theorem~\ref{thm:main}. The inclusion
$\mathbf{ID}(\rhd)\subset\IAA(\rhd)$ follows immediately as already mentioned. For the reverse inclusion,
let $\mu,\mu_n\in\mathbf P(\R)$ and let $k_n\to\infty$ be such that $\mu_n^{\rhd k_n}\Rightarrow\mu.$ It remains to show that $\mu\in\mathbf{ID}(\rhd)$.

For simplicity, we write $F_n:=F_{\mu_n}$, $F:=F_\mu$, and define
\[
V_n(z):=k_n(F_n(z)-z),
\quad z\in\C^+.
\] By the weak convergence assumption,
\[
F_n^{\circ k_n}
=
F_{\mu_n^{\rhd k_n}}
\longrightarrow F
\]
locally uniformly on $\C^+$.


\medskip
\noindent\textbf{Step~1: We have}
\[
\sup_{n\in\N}\frac{|V_n(iy)|}{y} \longrightarrow 0\quad \text{as}\quad y\longrightarrow+\infty. 
\]
Note that this implies $\{|V_n(iy)|\}_{n\in\N}$ is bounded for every sufficiently large $y > 0$; in particular $F_n(iy) \to iy$ as $n\to\infty$ and so by Lemma \ref{lem:Fconv} we have $\mu_n \Rightarrow \delta_0$.  

As the family $\{\mu_n^{\rhd k_n}: n\in \N\}$ is tight, by Lemma \ref{lem:Ftight} we have 
\begin{equation}\label{eq:tight}
\sup_{n \in \N}\left|\frac{F_n^{\circ k_n}(iy)-iy}{y}\right|
\longrightarrow 0 \quad \text{as} \quad y \longrightarrow +\infty. 
\end{equation}
In particular, there is $Y>0$  such that 
\[
\sup_{n \in \N,~y>Y}\left|\frac{\Im F_n^{\circ k_n}(iy)-y}{y}\right| < \log \frac{3}{2}. 
\]

We fix $n\in\N$ and $y>Y$ and define
\[
z_0:=iy,\quad
z_{j+1}:=F_n(z_j),\quad d_j:=z_{j+1}-z_j,\quad 0\le j\le k_n-1.
\]
Then
\[
F_n^{\circ k_n}(iy)-iy
=
\sum_{j=0}^{k_n-1}d_j.
\]

Suppose for the moment that $F_n-\mathrm{id}$ is nonconstant. By
Lemma \ref{lem:F},
\[
\Im(F_n(z)-z)\ge0,
\quad z\in\C^+.
\]
Since $\Im(F_n-\mathrm{id})$ is harmonic, the minimum principle gives
\[
\Im(F_n(z)-z)>0,
\quad z\in\C^+.
\]
Thus $F_n-\mathrm{id}$ is an analytic self-map of $\C^+$.
For $0\le j\le k_n-2$, applying
Lemma~\ref{lem:schwarz} to $f=F_n-\mathrm{id}$ at $z_j$ and
$z_{j+1}$ yields
\[
\frac{|d_{j+1}-d_j|^2}
{\Im d_{j+1}\Im d_j}
\le
\frac{|z_{j+1}-z_j|^2}
{\Im z_{j+1}\Im z_j}
=
\frac{|d_j|^2}
{\Im z_{j+1}\Im z_j}.
\]
Consequently,
\[\left|
\frac{d_{j+1}-d_j}{d_j}
\right|
\le
\sqrt{
\frac{\Im d_{j+1}\Im d_j}
{\Im z_{j+1}\Im z_j}
}
\le
\frac{\sqrt{\Im d_{j+1}\Im d_j}}{y}
\le
\frac{\Im d_{j+1}+\Im d_j}{2y}.\]
Here we used
\[
\Im z_{j+1}
=
\Im z_j+\Im d_j
>
\Im z_j
>\cdots>
\Im z_0=y,
\]
together with the arithmetic--geometric mean inequality.
Summing over $0\le j\le k_n-2$, we obtain
\[
\begin{aligned}
\sum_{j=0}^{k_n-2}
\left|
\frac{d_{j+1}-d_j}{d_j}
\right|\le
\frac{1}{2y}
\sum_{j=0}^{k_n-2}
(\Im d_{j+1}+\Im d_j)\le
\frac1y
\sum_{j=0}^{k_n-1}\Im d_j=
\frac{\Im F_n^{\circ k_n}(iy)-y}{y}<\log\frac{3}{2}.
\end{aligned}
\]
It follows that, for $0\le j\le k_n-1$,
\[
\begin{aligned}
\left|\frac{d_j}{d_0}-1\right|
&=
\left|
\prod_{\ell=0}^{j-1}\frac{d_{\ell+1}}{d_\ell}-1
\right|=
\left|
\prod_{\ell=0}^{j-1}
\left(
1+\frac{d_{\ell+1}-d_\ell}{d_\ell}
\right)-1
\right|\le
\prod_{\ell=0}^{j-1}
\left(
1+
\left|
\frac{d_{\ell+1}-d_\ell}{d_\ell}
\right|
\right)-1\\
&\le
\exp\left(
\sum_{\ell=0}^{j-1}
\left|
\frac{d_{\ell+1}-d_\ell}{d_\ell}
\right|
\right)-1\le
\frac12.
\end{aligned}
\]
 Then, for each $0\le j\le k_n-1$, we may write
\[
\frac{d_j}{d_0}=1+\eta_j,
\quad
|\eta_j|\le\frac12.
\]
Therefore,  

\begin{equation}\label{eq:V}
|F_n^{\circ k_n}(iy)-iy|
=
\left|\sum_{j=0}^{k_n-1}d_j\right|=
|d_0|
\left|k_n+
\sum_{j=0}^{k_n-1}\eta_j
\right|\ge
\frac12k_n|d_0| =\frac12|V_n(iy)|. 
\end{equation}
Note that the inequality $|F_n^{\circ k_n}(iy)-iy|\ge \frac12|V_n(iy)|$ trivially holds if $F_n-\text{id}$ is constant, as in this case $F_n^{\circ k_n}(iy)-iy = V_n(iy)$. Therefore inequality \eqref{eq:V} holds for all $n\in\N$ and $y>Y$. 
Step 1 follows by \eqref{eq:tight}. 

\medskip
\noindent\textbf{Step~2: The family $\{\mu_n^{\rhd m}: 0 \le m \le k_n,~n\in\N\}$ is tight.} 
Similar to \eqref{eq:V} we have
\begin{equation*}
|F_n^{\circ m}(iy)-iy|
=
\left|\sum_{j=0}^{m-1}d_j\right|=
|d_0|
\left|
\sum_{j=0}^{m-1}(1+\eta_j)
\right|\le
\frac32 m |d_0| \le\frac32|V_n(iy)|
\end{equation*}
for all $0\le m \le k_n$,  $n\in\N$ and $y>Y$. By Step 1 we have 
\[
\sup_{\substack{0\le m \le k_n \\ n\in\N}}\frac{|F_n^{\circ m}(iy)-iy|}{y} \longrightarrow 0 \quad \text{as}\quad y\longrightarrow +\infty. 
\]
Step 2 now follows by Lemma \ref{lem:Ftight}. 

\medskip
\noindent\textbf{Step~3: 
Construction of convolution roots.} We fix $m\in\N$. By Step 2, the family 
\[
\{\mu_n^{\rhd \lfloor \frac{k_n}{m} \rfloor}: n\in\N\}
\]
is tight. Therefore, it has a subsequence converging to a probability measure $\nu_m$. Since $k_n-m < m\lfloor \frac{k_n}{m} \rfloor \le k_n$, there is an integer $\ell \in \{0,1,\dots, m-1\}$ such that $m\lfloor \frac{k_n}{m}\rfloor -k_n = -\ell$ for infinitely many $n$'s within  the selected subsequence. For all such $n$'s we have 
\[
\mu_n^{\rhd m \lfloor \frac{k_n}{m} \rfloor} \rhd \mu_n^{\rhd \ell} = \mu_n^{\rhd k_n}.   
\]
The left side converges weakly to $\nu_m^{\rhd m}$ as $\mu_n\Rightarrow \delta_0$. Since the right side converges to $\mu$, we have $\mu = \nu_m^{\rhd m}$, i.e., $\mu \in \mathbf{ID}(\rhd)$.  

\begin{remark}
   The result in Step 2 under a slightly stronger assumption can be found in \cite[Lemma 3.2 (a)]{AW14}, which, however, contains a gap in its proof. More  specifically, $\Re F_{\rho_K}$ therein might go to infinity as $K\to+\infty$. Our proof of Step 2 not only generalizes \cite[Lemma 3.2 (a)]{AW14}, but also helps to fill the gap. 
\end{remark}

\section{General infinitesimal triangular arrays}\label{sec:proof2}

We prove Theorem~\ref{thm:ia}. First we consider
measures with finite second moment. Recall that a measure $\mu\in\mathbf P(\R)$ has finite second moment 
if and only if
\begin{equation}
F_\mu(z)
=
z-m_\mu-\int_\R\frac{1}{z-s}\d\tau_\mu(s),
\quad z\in\C^+,
\label{eq:finite-variance}
\end{equation}
for some $m_\mu\in\R$ and finite positive measure $\tau_\mu$. In this case, $m_\mu$ and $\tau_\mu(\R)$ coincide with the mean and variance of $\mu$, respectively. 
This follows from \cite[Proposition~2.2]{Maa92} by translation.

Suppose that $\mu \in \mathbf P(\R)$ has finite second moment with representation \eqref{eq:finite-variance}.  For $n\ge1$, define
\[
\nu_n:=\delta_{-(n+m_\mu)}\rhd\mu\rhd\delta_n.
\]
Then
\begin{align*}
F_{\nu_n}(z)=
F_{\delta_{-(n+m_\mu)}}\circ F_\mu\circ F_{\delta_n}(z)=
n+m_\mu+F_\mu(z-n)=
z-\int_{\R}\frac{1}{z-n-s}\d\tau_\mu(s).
\end{align*}
Since $\tau_\mu$ is a finite measure, the dominated convergence theorem gives $F_{\nu_n}(z)\to z$, and hence by Lemma \ref{lem:Fconv},  $\nu_n\Rightarrow\delta_0$.

We consider an infinitesimal triangular array with the $n$th row defined by 
\[
\underbrace{
\delta_{(n+m_\mu)/n^2},\ldots,\delta_{(n+m_\mu)/n^2}
}_{n^2\text{ copies}},
\quad
\nu_n,
\quad
\underbrace{
\delta_{-1/n},\ldots,\delta_{-1/n}
}_{n^2\text{ copies}}.
\]
Since
$
\nu_n\Rightarrow\delta_0,
$
this triangular array is infinitesimal. On the other hand,
\begin{align*}
&(F_{\delta_{(n+m_\mu)/n^2}})^{\circ n^2}
\circ F_{\nu_n}
\circ(F_{\delta_{-1/n}})^{\circ n^2}(z)
=-(n+m_\mu)+F_{\nu_n}(z+n)
=F_\mu(z).
\end{align*}
Thus the monotone convolution of every row is exactly $\mu$. Therefore,
every probability measure with finite second moment belongs to
$\mathbf{IA}(\rhd)$.

By \cite[Proposition~6.8(2)]{FHS20}, the class
$\mathbf{IA}(\rhd)$ is closed under weak convergence. Since compactly
supported probability measures are weakly dense in $\mathbf P(\R)$ and
have finite second moment, we conclude that
\[
\mathbf P(\R)\subset\mathbf{IA}(\rhd).
\]
The reverse inclusion is immediate from the definition, completing the proof.

\begin{remark}
\label{rem:exp}
Consider the symmetric Bernoulli measure
\[
\mu=\frac12(\delta_{-1}+\delta_1).
\]
Then
\[
G_\mu(z)=\frac{z}{z^2-1},
\qquad
F_\mu(z)=z-\frac1z,
\]
and $F_\mu'(i)=0$. Hence $F_\mu$ is not univalent on $\C^+$, whereas
Theorem~\ref{thm:ia} gives $\mu\in\mathbf{IA}(\rhd)$. Thus
\[
\mu\in\mathbf{IA}(\rhd)\setminus\mathbf{Univ}(\R),
\]
giving an explicit counterexample to Conjecture~6.10 of \cite{FHS20}.

This measure has $m_\mu=0$ and
$\tau_\mu=\delta_0$, so
\[
F_{\nu_n}(z)=z-\frac1{z-n}.
\]
This construction does not contradict
\cite[Theorem~6.9(2) and Remark~6.11]{FHS20}.
First,
\[
F_{\nu_n}(z+2\pi)
=
z+2\pi-\frac1{z+2\pi-n}
\neq F_{\nu_n}(z)+2\pi,
\]
so $\nu_n\notin\mathcal L$.
Second,
\[
F_{\nu_n}(z)
=
z-\int_{\R}\frac1{z-s}\d\delta_n(s),
\]
and hence $\operatorname{Var}(\nu_n)=\delta_n(\R)=1$. Since all the other entries
in the $n$th row are point masses, the maximal variance in that row is
$1$ for every $n$. Thus the hypothesis of
\cite[Theorem~6.9(2)]{FHS20} is not satisfied.
\end{remark}

\section{Multiplicative monotone convolution}\label{sec:unit}

We consider similar weak limit theorems concerning multiplicative monotone convolution $\circlearrowright$ on the unit circle. In fact, the proof is much simpler as the tightness  automatically holds due to the compactness of the unit circle. Recall that the multiplicative monotone convolution $\mu \circlearrowright \nu$ is the distribution of $UV$, where $U, V$ are unitary operators such that $U-I,V-I$ are monotonically independent and $U, V$ have distributions $\mu,\nu$, respectively. The reader is referred to \cite{Ber05,Fra06,FHS20} for further details and some  properties given below.

Let $\mathbf{P}(\mathbb{T})$ be the set of probability measures on the unit circle $\mathbb{T}$ and let $\mathbb D$ the unit disc of $\C$.   
For $\mu \in \mathbf{P}(\mathbb{T})$ let 
\[
\psi_\mu(z) := \int_\mathbb{T} \frac{z\zeta}{1-z\zeta} \d\mu(\zeta), \qquad \eta_\mu(z) := \frac{\psi_\mu(z)}{1+\psi_\mu(z)}, \qquad z \in \mathbb{D}.  
\]
We call $\eta_\mu$ the $\eta$-transform of $\mu$. 
The transform $\mu\mapsto \eta_\mu$ actually gives the  bijection  
\[
\mathbf{P}(\mathbb T)\longrightarrow\{\eta\colon \mathbb D\to \mathbb D~|~ \text{analytic}, \eta(0)=0\}. 
\]
It is known that multiplicative monotone convolution can be characterized by 
\[
\eta_{\mu \circlearrowright \nu} = \eta_\mu \circ \eta_\nu. 
\]
A measure $\mu\in\mathbf P(\mathbb{T})$ is called $\circlearrowright$-infinitely divisible if, for every $n\in\N$, there exists $\mu_n\in\mathbf P(\mathbb{T})$ such that
\[
\mu=\mu_n^{\circlearrowright n}.
\]
The class of all such measures is denoted by $\mathbf{ID}(\circlearrowright)$. In contrast to the case of additive monotone convolution $\rhd$, the $n$th roots $\mu_n$ may not be unique; for instance, $(\delta_{-1})^{\circlearrowright 2}=(\delta_{1})^{\circlearrowright 2}=\delta_1$. 

The weak convergence of a sequence $\{\mu_n\}_{n\in\N}\subset \mathbf{P}(\mathbb T)$ to a probability measure $\mu$ is denoted by $\mu_n\Rightarrow\mu.$
To discuss weak convergence,  the following fact is basic, see \cite[Lemma 2.11]{FHS20}.

\begin{lemma}\label{lem:conv_T}
    Let $\mu_n, \mu \in \mathbf P(\mathbb T), n\in\N$. Let $S\subset \mathbb D$ be a set containing an accumulation point in $\mathbb D$. The following are equivalent.
    \begin{enumerate}[label=\rm(\arabic*)]
    \item $ \mu_n\Rightarrow\mu$, 
    \item $ \eta_{\mu_n}\longrightarrow \eta_\mu$ locally uniformly on
$\mathbb D$, 
\item $\eta_{\mu_n}\longrightarrow \eta_\mu$ pointwise on $S$. 
    \end{enumerate}
\end{lemma}

Analogously to the case of additive monotone convolution, we consider the class $\IAA(\circlearrowright)$ of probability measures $\mu$ for which there exist a sequence $\{\mu_n\}_{n\in\N}\subset \mathbf P(\mathbb T)$ and a sequence $\{k_n\}_{n\in\N}\subset \N$ with $k_n \to \infty$ such that  
\[
\mu_n\Rightarrow \delta_1, \qquad \mu_n^{\circlearrowright k_n} \Rightarrow \mu. 
\]
Note that we do assume the infinitesimality $\mu_n \Rightarrow \delta_1$ unlike the additive case. This is because the infinitesimality does not hold automatically; for example, $(\delta_{-1})^{\circlearrowright 2n}=\delta_1$ and $\haar^{\circlearrowright n} =\haar$, where $\haar$ is the normalized Haar measure on $\mathbb T$.

The following fact is a consequence of \cite[Theorem 4.7]{Ber05}. Since a detailed proof is omitted in \cite[Theorem 4.7]{Ber05}, we give a self-contained proof  for the reader's convenience.  
 
\begin{lemma}\label{lem:ID}
Let $\mu\in \mathbf{ID}(\circlearrowright)\setminus\{\haar\}$. Then there exists a sequence $\{\mu_n\}_{n\in \N}\subset\mathbf{P}(\mathbb T)$  such that $\mu=\mu_n^{\circlearrowright n}$ for all $n$ and $\mu_n\Rightarrow\delta_1$. 
\end{lemma}

\begin{proof}

We set $f:=\eta_\mu$ and $a:=f'(0)$. 

\medskip
\noindent\textbf{Step 1: $0<|a|\le 1$.} Suppose that $a=0$, and let
$d\ge2$ be the order of the zero of $f$ at $0$. Choose $N\in\N$
such that $2^N>d$. Since $\mu$ is $\circlearrowright$-ID, there exists an analytic map $g\colon\mathbb{D}\to\mathbb{D}$
such that $g(0)=0$ and $g^{\circ N}=f$. Then $g\not\equiv0$ and
$(g'(0))^N=f'(0)=0$. Thus the order $\tilde{d}$ of the zero of $g$
at $0$ satisfies $\tilde{d}\ge2$. Since orders of zeros multiply
under composition, we obtain
\[
d=\tilde{d}^{N}\ge2^N>d,
\]
a contradiction. Hence $|a|>0$. The other bound $|a|\le1$ follows from the Schwarz lemma, as $f\colon\mathbb D\to \mathbb D$ and $f(0)=0$. 

\medskip
\noindent\textbf{Step 2: The linear case.}
If $f(z)=az$, fix $\theta\in(-\pi, \pi]$ with $a=|a| e^{i\theta}$.
For each $n\ge1$, the map
\[
z\longmapsto |a|^{1/n}e^{i\theta/n}z
\]
is an analytic self-map of $\mathbb{D}$ fixing $0$  since $|a|\le1$. Taking the corresponding measures $\mu_n$, which are actually Poisson kernels or delta measures, proves the assertion in the linear case.

\medskip
\noindent\textbf{Step 3: Comparing the first nonlinear coefficients.}
Assume now that $f$ is not linear. Schwarz's lemma and its equality case imply $0<|a|<1$. Write
\[
f(z)=az+Az^{q+1}+O(z^{q+2}),
\quad A\ne0,
\quad q\ge1.
\]
By the assumption of the infinite divisibility, for every $n\ge1$ there is an analytic self-map $g_n$ of $\mathbb{D}$, fixing $0$, such that
\[
g_n^{\circ(qn)}=f. 
\]
Let $b_n:=g_n'(0).$ Then
\[
(b_n)^{qn}=a,
\quad |b_n|=|a|^{1/(qn)}\in(0,1).
\]
Since $g_n\circ f=f\circ g_n$ and $f(z)=az+O(z^{q+1})$,
comparison of the coefficients of $z^2,\ldots,z^q$, together with $a^k\ne a$ for $2\le k\le q$, gives
\[
g_n(z)=b_nz+B_nz^{q+1}+O(z^{q+2}).
\]
At degree $q+1$, the two compositions are
\begin{align*}
g_n(f(z))&=ab_nz+(Ab_n+B_na^{q+1})z^{q+1}+O(z^{q+2}),\\
f(g_n(z))&=ab_nz+(aB_n+Ab_n^{q+1})z^{q+1}+O(z^{q+2}).
\end{align*}
Equating these coefficients and rearranging yields
\begin{equation}\label{coef}
a(1-a^q)B_n=Ab_n(1-b_n^q).
\end{equation}

\medskip
\noindent\textbf{Step 4: Estimating coefficients.}
Define $h_n(z)=g_n(z)/z$ for $z\ne0$, and $h_n(0)=b_n$. This extension is analytic. Schwarz's lemma gives $|h_n(z)|\le1$.
Since $|h_n(0)|=|b_n|<1$, the maximum modulus principle gives $|h_n(z)|<1$ throughout $\mathbb{D}$, i.e., $h_n$ is a self-map of $\mathbb D$. We further compose $h_n$ and the automorphism $z\mapsto (z-b_n)/(1-\overline{b_n}z)$ of $\mathbb D$:  
\[
k_n(z):=\frac{h_n(z)-b_n}{1-\overline{b_n}h_n(z)} = \frac{B_n}{1-|b_n|^2}z^q+O(z^{q+1}). 
\]
The repeated application of Schwarz's lemma to $k_n(z), k_n(z)/z,\dots, k_n(z)/z^{q-1}$ yields that these are all self-maps of $\mathbb D$ and $|k_n(z)/z^{q-1}|\le |z|$, which further implies
\begin{equation}\label{schwarz}
|B_n|\le1-|b_n|^2.
\end{equation}
Combining \eqref{coef} and \eqref{schwarz}, and using $|b_n|= |a|^{1/(qn)}\ge|a|$, we obtain
\[
|1-b_n^q|
\le\frac{|a||1-a^q|}{|A||b_n|}(1-|b_n|^2)
\le\frac{|1-a^q|}{|A|}(1-|b_n|^2).
\]
As $|b_n|=|a|^{1/(qn)}\to1$ we obtain $(g_n^{\circ q})'(0)=b_n^q\to 1$. This implies $g_n^{\circ q}\to \mathrm{id}$ by a compactness argument and the extremal case of Schwarz's lemma. 
Let $\mu_n$ be the probability measure with $\eta_{\mu_n}=g_{n}^{\circ q}$.
Since $(g_{n}^{\circ q})^{\circ n}=g_{n}^{\circ(qn)}=f=\eta_\mu$, the transform identities and uniqueness give
$
\mu_n^{\circlearrowright n}=\mu
$. 
Finally, by Lemma \ref{lem:conv_T}, $\eta_{\mu_n}=g_{n}^{\circ q}\to\mathrm{id}$ implies $\mu_n\Rightarrow\delta_1$ as required.
\end{proof}

\begin{theorem}\label{thm:unit}
    We have 
    \[
    \IAA(\circlearrowright) = \mathbf{ID}(\circlearrowright). 
    \]
\end{theorem}
\begin{proof}
 For the inclusion $\IAA(\circlearrowright) \subset \mathbf{ID}(\circlearrowright)$, we can simply follow the lines of Step 3 in Section \ref{sec:proof}. Note that by Helly's selection theorem (or more general Prokhorov's theorem), any family of probability measures on $\mathbb T$ is tight.  

 As for the converse inclusion,  Lemma \ref{lem:ID} yields $\mathbf{ID}(\circlearrowright)\setminus\{\haar\} \subset \IAA(\circlearrowright)$. The Haar measure $\haar$ can be approximated by Poisson kernels $\mu_n$ defined by $\eta_{\mu_n}(z)= a_n z$ with $a_n\in \C,~|a_n|<1$, such that $a_n\to1$ and  $(a_n)^n\to0$, e.g., $a_n=1-n^{-1/2}$. Then $\eta_{\mu_n}(z)\to z =\eta_{\delta_1}(z)$ and $\eta_{\mu_n^{\circlearrowright n}}(z)=(a_n)^n z \to 0 = \eta_\haar(z)$; therefore we obtain $\haar \in \IAA(\circlearrowright)$. 
\end{proof}

\begin{remark}\label{rem:unit}
    Combining this theorem with the known inclusions in \cite[Section 7]{FHS20}, we obtain
\[
    \IAA(\circlearrowright)
    =\mathbf{ID}(\circlearrowright)
    \subsetneq\mathbf{Univ}(\mathbb T)\cup\{\haar\}
    =\mathbf{IA}(\circlearrowright)
    \subsetneq\mathbf P(\mathbb T), 
\]
where $\mathbf{Univ}(\mathbb T)$ is the set of probability measures with univalent $\eta$-transform and $\mathbf{IA}(\circlearrowright)$ is the collection of all $\mu \in \mathbf{P}(\mathbb T)$ such that there exists a family $\{\mu_{n,j}: 1\le j \le k_n, ~n\in\N\}$ in $\mathbf{P}(\mathbb T)$ with $k_n\in\N$, $k_n\to\infty$ and 
\begin{align*}
&\lim_{n\to\infty}\max_{1\le j \le k_n}\mu_{n,j}(\{e^{i\theta}:\epsilon<\theta <2\pi -\epsilon\}) =0, \qquad \forall \epsilon \in (0,\pi), \qquad \text{(infinitesimality)}\\
&\mu_{n,1}\circlearrowright \mu_{n,2}\circlearrowright \cdots \circlearrowright \mu_{n,k_n} \Rightarrow \mu. 
\end{align*}
These inclusions are not completely parallel to  the additive case  \eqref{eq:inclusions}. Note that a proof of the strictness of the inclusion $\mathbf{ID}(\circlearrowright)
    \subsetneq\mathbf{Univ}(\mathbb T)\cup\{\haar\}$ is not provided in the literature, but the geometric idea in \cite[Remark 6.15]{FHS20} works to construct an example of measures  in $\mathbf{Univ}(\mathbb T)\setminus \mathbf{ID}(\circlearrowright)$. 
\end{remark}

\section*{Statement on AI use}
The manuscript is written with the help of GPT-6 Astra Ultra (OpenAI) for exploratory
mathematical discussions. In particular, the AI suggested considering
the Schwarz--Pick lemma in developing certain estimates and proposed
a special case of a potential counterexample to the conjecture $\mathbf{IA}(\rhd)=\mathbf{Univ}(\R)$. The authors
substantially refined the resulting arguments and generalized
the counterexample. They independently verified all mathematical
results and take full responsibility for the contents of this paper.

\end{document}